\documentclass[11pt,reqno]{amsart}

\usepackage{latexsym,amsmath,amssymb}
\begin{document}

\title{Supercommutator algebras of right (Hom-)alternative superalgebras}


\author{A. Nourou Issa} 
         




\begin{abstract}
The supercommutator algebra of a right alternative superalgebra is a Bol superalgebra. Hom-Bol superalgebras are defined
and it is shown that they are closed under even self-morphisms. Any Bol superalgebra along with any even self-morphism
is twisted into a Hom-Bol superalgebra. The supercommutator algebra of a right Hom-alternative superalgebra has a natural
Hom-Bol superalgebra structure. 
\end{abstract}

\date{}

\maketitle
\section{Introduction}
A {\it right alternative algebra} is an algebra satisfying the {\it right alternative identity} $(xy)y = x(yy)$. These algebras
were first considered in \cite{Alb}. For further studies on right alternative algebras, one may refer to \cite{Kl}, \cite{Sk}, 
\cite{Thed}.
\par
It turns out that right alternative algebras have close relations with a type of binary-ternary algebras called {\it Bol algebras}
which were introduced in \cite{Mik1} (see also references therein). In fact, it is proved \cite{Mik2} that any right alternative
algebra has a natural Bol algebra structure.
\par
The general theory of superalgebras started with the introduction of $\mathbb{Z}_2$-graded Lie algebras (i.e. Lie superalgebras)
coming from physics (see \cite{Kac1}, \cite{Sch} and references therein for basics on Lie superalgebras). The $\mathbb{Z}_2$-graded
generalization of algebras is first extended to Jordan algebras in \cite{Kac2}. Next, alternative superalgebras were introduced in
\cite{ZS} whereas Maltsev superalgebras were introduced in \cite{Sh1}. A $\mathbb{Z}_2$-graded generalization of Bol algebras is 
considered in \cite{Ruk}.
\par
With the introduction of Hom-Lie algebras (see \cite{HLS}, \cite{LS1}, \cite{LS2}) began studies of Hom-type generalizations of usual
algebras. Apart from Hom-Lie algebras, first Hom-type algebras were defined in \cite{MS} while Hom-alternative and Hom-Jordan
algebras were defined in \cite{Ma} (see also \cite{Yau3} where Maltsev algebras were defined) and the Hom-type generalization of right
alternative algebras was considered in \cite{Yau2}. It should be observed that, in general, the twisting map in a Hom-algebra is neither
injective nor surjective and when the twisting map is the identity map, then one recovers the ordinary (untwisted) algebraic structure. 
So ordinary algebras are viewed as Hom-algebras with the identity map as twisting map. Moving further in the theory of Hom-algebras, 
the twisting principle of algebras
is extended to binary-ternary algebras in \cite{Iss1} and next Hom-Bol algebras were defined in \cite{AI1}. As in the case of right
alternative algebras, it is shown in \cite{AI2} that a Hom-Bol algebra structure can be defined on any multiplicative right Hom-alternative 
algebra.
\par
In this paper we extend, but with a different approach, the result in \cite{Mik2} to the cases of right alternative superalgebras and 
right Hom-alternative superalgebras.
Some basics on superalgebras are reminded in section 2. In section 3 it is proved that any right alternative superalgebra has a natural
Bol superalgebra structure. In section 4, in order to deal with the Hom-version of the results from section 3, we first define Hom-Bol
superalgebras and next, in section 5, prove that right Hom-alternative superalgebras are in fact Hom-Bol superalgebras.
\par
All vector spaces and algebras are considered over a fixed ground field of characteristic not $2$ or $3$.
\section{Preliminaries}
A {\it superspace} (or a $\mathbb{Z}_2$-{\it graded space}) $V$ is  a direct sum $V = V_0\oplus V_1$, where $V_i$ are vector spaces. An element 
$x \in V_i$ ($i \in \mathbb{Z}_2$) is said to be {\it homogeneous of degree} $i$ and the degree of $x$ will be denoted by $\bar{x}$. \\
\par
{\bf Definition 2.1}
 (i) Let $f: A \rightarrow A'$ be a linear map, where $A = A_0 \oplus A_1$ and $A' = A'_0 \oplus A'_1$ are superspaces. The
map $f$ is said to be {\it even} (resp. {\it odd}) if $f(A_i) \subset A'_i$ (resp. $f (A_i ) \subset A'_{i+1})$ for $i = 0,1$. 
\par
(ii) A ({\it multiplicative}) $n$-{\it ary Hom-superalgebra} is a triple ($A, \{\cdot , \cdots , \cdot \}, \alpha$) consisting of a 
superspace $A = A_0 \oplus A_1$, an $n$-linear map $\{\cdot , \cdots , \cdot \}: A^{\otimes n} \rightarrow A$ such that 
$\{A_i , \cdots , A_s \} \subset A_{i+ \cdots +s}$, and an even linear map $\alpha: A \rightarrow A$ such that
$\alpha (\{ x_1 , \cdots , x_n \}) = \{ \alpha (x_1), \cdots , \alpha (x_n) \}$ (multiplicativity).\\
\par
One observes that if $\alpha = Id$ (the identity map), we get the corresponding definition of an $n$-{\it ary superalgebra}. One also notes
that $\overline{\alpha (x)} = \bar{x}$ for all homogeneous $x \in A$.
\par
We will be interested in binary ($n=2$), ternary ($n=3$) and binary-ternary Hom-algebras (i.e. Hom-algebras with binary and ternary
operations). For convenience, throughout this paper we assume that all Hom-(super)algebras are multiplicative.\\
\par
{\bf Definition 2.2} (\cite{AM}). Let ${\mathcal A} := (A, *, \alpha)$ be a binary Hom-superalgebra. The {\it supercommutator Hom-algebra}
(or the {\it minus Hom-superalgebra}) of $\mathcal A$ is the Hom-superalgebra ${\mathcal A}^{-} := (A, [\cdot , \cdot], \alpha)$, 
where $[x,y] := {1 \over 2}(x * y - (-1)^{\bar{x} \bar{y}} y * x)$ for all homogeneous $ x,y \in A$. The product $[\cdot , \cdot]$ is called the 
{\it supercommutator bracket}. \\
\par
In the sequel, we will also denote by ``$[\cdot , \cdot]$`` the binary operation in (Hom-)superalgebras. Besides the supercommutator of elements
in Hom-superalgebras, one also considers the {\it super-Jordan product}
$$ x \circ y := {1 \over 2}(xy + (-1)^{\bar{x} \bar{y}} yx),$$
the {\it Hom-Jordan associator}
$$as^{J}_{\alpha}(x,y,z) := xy \cdot \alpha (z) + \alpha (x) \cdot yz $$
and the {\it Hom-associator} $as_{\alpha}(x,y,z)$ defined as
$$as_{\alpha}(x,y,z) := xy \cdot \alpha (z) - \alpha (x) \cdot yz $$
for $x,y,z$ in the given Hom-superalgebra \cite{MS}. When $\alpha = Id$ we recover the usual {\it Jordan associator} $as^{J}(x,y,z)$ and 
{\it associator} $as(x,y,z)$ respectively in usual algebras. The Hom-superalgebra ${\mathcal A}^{+} := (A, \circ , \alpha)$ is usually called
the {\it plus Hom-superalgebra} of ${\mathcal A} := (A, *, \alpha)$.\\
\par
{\bf Definition 2.3} A Hom-superalgebra $A$ is said to be {\it right Hom-alternative} if \\
\par
$as_{\alpha}(x,y,z) = - (-1)^{\bar{y} \bar{z}} as_{\alpha}(x,z,y)$ ({\it right superalternativity}) \hfill (2.1)\\
\\
for all $x,y,z \in A$. \\
\par
Likewise is defined a left Hom-alternative superalgebra. A Hom-superalgebra that is both right and left Hom-alternative is said to be 
{\it Hom-alternative}. Of course, for $\alpha = Id$, one gets the definition of an alternative superalgebra given in \cite{ZS}.
\par
From (2.1), expanding Hom-associators, it is easily seen that (2.1) is equivalent to \\
\par
$\alpha (x)(yz + (-1)^{\bar{y} \bar{z}} zy) = (xy)\alpha (z) + (-1)^{\bar{y} \bar{z}} (xz)\alpha (y)$. \hfill (2.2) 
\section{Right alternative superalgebras and Bol superalgebras}
In this section we prove that the commutator algebra of a right alternative superalgebra is a Bol superalgebra. First we recall the following\\
\par
{\bf Definition 3.1} (\cite{Ruk}). A {\it Bol superalgebra} is a triple $(A, [\cdot , \cdot], \{ \cdot, \cdot, \cdot \})$ in which $A$ is a
superspace, $[\cdot , \cdot]$ and $\{ \cdot, \cdot, \cdot \}$ are binary and ternary operations on $A$ such that \\
\par
(SB1) $[x,y] = - (-1)^{\bar{x} \bar{y}} [y,x]$,
\par
(SB2) $\{x,y,z\} = - (-1)^{\bar{x} \bar{y}} \{y,x,z\}$,
\par
(SB3) $\{x,y,z\} + (-1)^{\bar{x} (\bar{y} + \bar{z})} \{y,z,x\} + (-1)^{\bar{z} (\bar{x} + \bar{y})}\{z,x,y\} = 0$,
\par
(SB4) $\{x,y,[u,v]\} = [\{x,y,u\},v] + (-1)^{\bar{u} (\bar{x} + \bar{y})}[u,\{x,y,v\}]$
\par
\hspace{3.5cm}$+ (-1)^{(\bar{x} + \bar{y}) (\bar{u} + \bar{v})} (\{u,v,[x,y]\} - [[u,v],[x,y]])$,
\par
(SB5) $\{x,y,\{u,v,w\} \} = \{ \{x,y,u\}, v,w\} + (-1)^{\bar{u} (\bar{x} + \bar{y})} \{u, \{x,y,v\}, w\}$
\par
\hspace{4.2cm}$+ (-1)^{(\bar{x} + \bar{y}) (\bar{u} + \bar{v})} \{u,v, \{x,y,w\} \}$ \\
\\
for all homogeneous $u,v,w,x,y,z \in A$.\\
\par
Clearly, any Bol superalgebra with zero odd part is a (left) Bol algebra. If $[x,y] = 0$ for all homogeneous $x,y \in A$, then 
$(A, [\cdot , \cdot], \{ \cdot, \cdot, \cdot \})$ reduces to a {\it Lie supertriple system} $(A,\{ \cdot, \cdot, \cdot \})$.\\
\par
{\it Example 3.1} Let $A = A_0 \oplus A_1$ be a superspace where $A_0$ is a $2$-dimensional vector space with basis $\{i,j\}$ and  $A_1$ a 
$1$-dimensional vector space with basis $\{k\}$. Define on $A$ the following binary and ternary nonzero products:
\par
$[i,j]=j$, $[i,k]=k$,
\par
$[j,i]=-j$,
\par
$[k,i]=k$, $[k,k]=j$,
\par
$\{i,j,i\}=-j$, $\{i,k,i\}=-k$,
\par
$\{j,i,i\}=j$,
\par
$\{k,i,i\}=k$. \\
Then it could be checked that $(A, [\cdot , \cdot], \{ \cdot, \cdot, \cdot \})$, with the multiplication table as above, is a ($3$-dimensional) 
Bol superalgebra. Observe that $(A, [\cdot , \cdot])$ is a $3$-dimensional Maltsev superalgebra \cite{AE} and the table for the ternary product 
as above is obtained using the $\mathbb{Z}_2$-graded version of the ternary product that produces a Bol algebra from a Maltsev algebra \cite{Mik1}.
\par
In \cite{Mik2} it is proved that on any right alternative algebra one may define a Bol algebra structure. The $\mathbb{Z}_{2}$-graded version of this 
result is given by the following \\
\par
{\bf Theorem 3.1} {\it The supercommutator algebra of any right alternative superalgebra is a Bol superalgebra}. \\
\par
{\it Proof} Let $A$ be a right alternative superalgebra. Then $(A, \circ )$ is a Jordan superalgebra \cite{Sh2}. Now define on $(A, \circ )$
a ternary product
\par
$[x,y,z] := 2(-1)^{\bar{x}(\bar{y} + \bar{z})} as^{J}(y,z,x)$. \\
Then one checks that $(A, [\cdot, \cdot, \cdot])$ is a Lie supertriple system. So is $(A, \{ \cdot, \cdot, \cdot \})$, where $\{x,y,z\} = 
(-1)^{\bar{x}(\bar{y} + \bar{z})} as^{J}(y,z,x)$. Now, using specific properties of right alternative superalgebras, one gets that
$(A, [\cdot, \cdot], \{ \cdot, \cdot, \cdot \})$ is a Bol superalgebra, where $[\cdot, \cdot]$ is the supercommutator operation on 
$A$. \hfill $\square$ 
\section{Hom-Bol superalgebras. Construction theorems and example}
In \cite{AI1} Hom-Bol algebras are defined. In this section we define Hom-Bol superalgebras as a generalization both of Bol superalgebras \cite{Ruk}
and Hom-Bol algebras \cite{AI1}. Next we point out  some construction theorems.\\
\par
{\bf Definition 4.1} A {\it Hom-Bol superalgebra} is a quadruple ${\mathcal{A}}_{\alpha} := (A, [\cdot , \cdot], \{ \cdot, \cdot, \cdot \}, \alpha)$
where $A$ is a superspace, $[\cdot , \cdot]$ (resp. $\{ \cdot, \cdot, \cdot \}$) is a binary (resp. ternary) operation on $A$ such that\\
\par
(SHB1) $\alpha ([x,y]) = [\alpha (x), \alpha (y)]$,
\par
(SBH2) $\alpha (\{x,y,z\}) = \{ \alpha (x), \alpha (y), \alpha (z) \}$,
\par
(SBH3) $[x,y] = - (-1)^{\bar{x} \bar{y}} [y,x]$,
\par
(SBH4) $\{x,y,z\} = - (-1)^{\bar{x} \bar{y}} \{y,x,z\}$,
\par
(SBH5) $\{x,y,z\} + (-1)^{\bar{x} (\bar{y} + \bar{z})} \{y,z,x\} + (-1)^{\bar{z} (\bar{x} + \bar{y})}\{z,x,y\} = 0$,
\par
(SBH6) $\{ \alpha (x), \alpha (y), [u,v] \} = [\{x,y,u\},{{\alpha}^{2}}(v)] + (-1)^{\bar{u} (\bar{x} + \bar{y})}[{{\alpha}^{2}}(u),\{x,y,v\}]$
\par
\hspace{1.5cm}$+ (-1)^{(\bar{x} + \bar{y}) (\bar{u} + \bar{v})} (\{\alpha(u),\alpha(v),[x,y]\} - [[\alpha(u),\alpha(v)],[\alpha(x),\alpha(y)]])$,
\par
(SBH7) $\{{{\alpha}^{2}}(x),{{\alpha}^{2}}(y),\{u,v,w\} \} = \{ \{x,y,u\}, {{\alpha}^{2}}(v),{{\alpha}^{2}}(w)\}$
\par
\hspace{5.9cm}$+ (-1)^{\bar{u} (\bar{x} + \bar{y})} \{{{\alpha}^{2}}(u), \{x,y,v\}, {{\alpha}^{2}}(w)\}$
\par
\hspace{5.9cm}$+ (-1)^{(\bar{x} + \bar{y}) (\bar{u} + \bar{v})} \{{{\alpha}^{2}}(u),{{\alpha}^{2}}(v), \{x,y,w\} \}$ \\
\\
for all homogeneous $u,v,w,x,y,z \in A$.\\
\par
We observe that for $\alpha = Id$, any Hom-Bol superalgebra reduces to a Bol superalgebra and a Hom-Bol superalgebra with a zero odd part is a 
Hom-Bol algebra. If $[x,y]=0$ for all homogeneous $x,y \in A$, one gets a {\it Hom-Lie supertriple system} $(A,\{ \cdot, \cdot, \cdot \}, {\alpha}^2)$. \\
\par
{\bf Theorem 4.1} {\it Let ${\mathcal{A}}_{\alpha} := (A, [\cdot , \cdot], \{ \cdot, \cdot, \cdot \}, \alpha)$ be a Hom-Bol superalgebra and 
$\beta : A \rightarrow A$ an even self-morphism of ${\mathcal{A}}_{\alpha}$ such that $\alpha \circ \beta = \beta \circ \alpha$. Let ${\beta}^0 = Id$
and ${\beta}^n = \beta \circ {\beta}^{n-1}$ for any integer $n \geq 0$ and define on $A$ a binary operation $[\cdot, \cdot]_{{\beta}^n}$ and a
ternary operation $\{ \cdot, \cdot, \cdot \}_{{\beta}^n}$ by
\par
$[x, y]_{{\beta}^n} := {\beta}^n ([x,y])$,
\par
$\{ x,y,z \}_{{\beta}^n} := {\beta}^{2n} (\{x,y,z\})$. \\
Then ${\mathcal{A}}_{{\beta}^n} := (A, [\cdot, \cdot]_{{\beta}^n}, \{ \cdot, \cdot, \cdot \}_{{\beta}^n}, {\beta}^n \circ \alpha)$ is a Hom-Bol
superalgebra}. \\
\par
{\it Proof} The proof is similar to that of Theorem 3.2 in \cite{AI1}. \hfill $\square$ \\ 
\par
From Theorem 4.1 we get the following extension of the Yau's twisting principle \cite{Yau1} giving a construction of Hom-Bol superalgebras from
Bol superalgebras.\\
\par
{\bf Corollary 4.1} {\it Let $(A, [\cdot , \cdot], \{ \cdot, \cdot, \cdot \})$ be a Bol superalgebra and $\beta$ an even self-morphism of
$(A, [\cdot , \cdot], \{ \cdot, \cdot, \cdot \})$. Define on $A$ a binary operation $[\cdot, \cdot]_{\beta}$ and a ternary operation 
$\{ \cdot, \cdot, \cdot \}_{\beta}$ by
\par
$[x, y]_{\beta} := {\beta} ([x,y])$,
\par
$\{ x,y,z \}_{\beta} := {\beta}^{2} (\{x,y,z\})$. \\
Then ${\mathcal{A}}_{\beta} := (A, [\cdot, \cdot]_{\beta}, \{ \cdot, \cdot, \cdot \}_{\beta}, \beta )$ is a Hom-Bol superalgebra. Moreover, if 
$(A', [\cdot , \cdot]', \{ \cdot, \cdot, \cdot \}')$ is another Bol superalgebra, ${\beta}'$ an even self-morphism of $(A', [\cdot , \cdot]', \{ \cdot, \cdot, \cdot \}')$
and if $f:A \rightarrow A'$ is a Bol superalgebra even morphism satisfying $f \circ \beta = {\beta}' \circ f$, then 
$f:{\mathcal A}_{\beta} \rightarrow {\mathcal{A}'}_{{\beta}'}$ is a morphism of Hom-Bol superalgebras, where 
${\mathcal{A}'}_{{\beta}'} := (A', [\cdot , \cdot{]'}_{{\beta}'}, \{ \cdot, \cdot, \cdot {\}'}_{{\beta}'}, {\beta}')$}.\\
\par
{\it Proof} The first part of the corollary comes from Theorem 4.1 when $n=1$. The second part is proved in a similar way as Corollary 4.5 
in \cite{Iss1}. \hfill $\square$ \\ 
\par 
{\it Example 4.1} Let  $A = A_0 \oplus A_1$ be a superspace over a field of characteristic not $2$ where $A_0$ is a $1$-dimensional vector space 
with basis $\{i\}$ and  $A_1$ a $2$-dimensional vector space with basis $\{j,k\}$. Define on $A$ the following only nonzero products on basis elements:
\par
$i*j=k$;
\par
$j*i=k$, $j*k=2i$;
\par
$k*j=4i$.\\
Then $(A, *)$ is a right alternative superalgebra \cite{Sh2}. Now consider on $(A, *)$ the supercommutator $[\cdot , \cdot]$ and the ternary operation
defined as
\par
$\{x,y,z\} := (-1)^{\bar{x}(\bar{y} + \bar{z})} as^{J}(y,z,x)$.\\
Then it could be checked that $(A,[\cdot,\cdot],\{ \cdot,\cdot,\cdot \})$ is a Bol superalgebra, where the only nonzero products are:
\par
$[j,k] = 6i$;
\par
$[k,j] = 6i$;
\par
$\{i,j,j\} = 4i$;
\par
$\{j,i,j\} = -4i$, $\{j,j,i\} = -8i$, $\{j,j,k\} = -8k$, $\{j,k,j\} = 4k$;
\par
$\{k,j,j\} = 4k$.\\
Next define a linear map $\beta : A \rightarrow A$ by setting
\par
$\beta (i) = ai$, $\beta (j) = j+bk$, $\beta (k) = ak$\\
with $a \neq 0$. Then it is easily seen that $\beta$ is an even self-morphism of $(A,[\cdot,\cdot],\{ \cdot,\cdot,\cdot \})$ and Corollary 4.1 implies
that ${\mathcal{A}}_{\beta} := (A, [\cdot, \cdot]_{\beta}, \{ \cdot, \cdot, \cdot \}_{\beta}, \beta )$ is a Hom-Bol superalgebra with products given as
\par
$[j,k]_{\beta} = 6ai$;
\par
$[k,j]_{\beta} = 6ai$;
\par
$\{i,j,j\}_{\beta} = 4a^{2}i$;
\par
$\{j,i,j\}_{\beta} = -4a^{2}i$, $\{j,j,i\}_{\beta} = -8a^{2}i$, $\{j,j,k\}_{\beta} = -8a^{2}k$, $\{j,k,j\}_{\beta} = 4a^{2}k$;
\par
$\{k,j,j\}_{\beta} = 4a^{2}k$.\\
\par
The notion of an $n$th derived (binary) Hom-algebra of a given Hom-algebra is first introduced in \cite{Yau3} and the closure of a given type of 
Hom-algebras under taking $n$th derived
Hom-algebras is a property that is characteristic of the variety of Hom-algebras. Later on, this notion is extended to binary-ternary Hom-algebras \cite{AI1} 
(for binary-ternary Hom-superalgebras \cite{GI}, the notion is the same as in the case of binary-ternary Hom-algebras).\\
\par
{\bf Definition 4.2}  (\cite{GI}) Let ${\mathcal A} := (A,*,\{ \cdot,\cdot,\cdot \},\alpha)$ be a binary-ternary Hom-superalgebra and $n\geq 0$ an integer. 
Define on $A$ the $n$th derived 
binary operation $*^{(n)}$ and the $n$th derived ternary operation $\{ \cdot,\cdot,\cdot \}^{(n)}$ by
\par
$x *^{(n)} y := \alpha^{{2^n}-1}(x * y)$,
\par
$\{x,y,z\}^{(n)} := \alpha^{{2^{n+1}}-2}(\{x,y,z\})$, \\
for all homogeneous $x,y,z$ in A. Then ${\mathcal A}^{(n)} := (A, *^{(n)},\{ \cdot,\cdot,\cdot \}^{(n)}, \alpha^{2^{n}})$ is called the $n$th 
{\it derived  (binary-ternary) 
Hom-superalgebra} of $\mathcal A$. \\
\par 
As for Hom-Bol algebras, the category of Hom-Bol superalgebras is closed under taking derived Hom-superalgebras as stated in the following\\
\par
{\bf Theorem 4.2} {\it Let ${\mathcal A} := (A,[\cdot,\cdot],\{ \cdot,\cdot,\cdot \},\alpha)$ be a Hom-Bol superalgebra.
Then, for each $n \geq 0$, the $n$th derived Hom-superalgebra ${\mathcal A}^{(n)}$ is a Hom-Bol superalgebra. In
particular, the $n$th derived Hom-superalgebra of a Hom-Lie supertriple system is a Hom-Lie supertriple system}. \\
\par
{\it Proof} The proof of the first part of the theorem is similar to that of Theorem 3.5 in \cite{AI1}, and the second part follows immediately.
\hfill $\square$  
\section{Hom-Bol superalgebra structures on right Hom-alternative superalgebras}
In this section we prove a $\mathbb{Z}_2$-graded generalization of results connecting right Hom-alternative algebras and Hom-Bol algebras \cite{AI2}.\\
\par
{\bf Lemma 5.1} {\it Let ${\mathcal A} := (A,*, \alpha)$ be a multiplicative right Hom-alternative superalgebra. Then ${\mathcal A}$ is 
Hom-Jordan-admissible, i.e. ${\mathcal A}^{+}$ is a Hom-Jordan superalgebra}.\\
\par
{\it Proof} A proof comes from a generalization to right alternative case of the proof of Theorem 6.1 in \cite{AAM}. \hfill $\square$ \\ 
\par
{\bf Lemma 5.2} {\it Let $(A,\circ, \alpha)$ be a multiplicative Hom-Jordan superalgebra. If define on $A$ a ternary product as}
\par
$[x,y,z] := 2 (-1)^{\bar{x}(\bar{y} + \bar{z})} as^{J}_{\alpha}(y,z,x)$ \hfill (5.1) \\
{\it for all homogeneous $x,y,z \in A$, then $(A, [ \cdot, \cdot, \cdot ], {\alpha}^2)$ is a Hom-Lie supertriple system}.\\
\par
{\it Proof} If define on $(A,\circ, \alpha)$ a ternary product
\par
$(x,y,z) := (x \circ y) \circ \alpha (z) + (-1)^{\bar{x}\bar{y} + \bar{z}(\bar{x} + \bar{y})} (z \circ y) \circ \alpha (x) - (-1)^{\bar{x}\bar{y}} 
\alpha (y) \circ (x \circ z)$,\\
then it is verified that $(A, (\cdot, \cdot, \cdot), {\alpha}^2)$ is a Hom-Jordan supertriple system. Next, defining on $A$ the ternary product
\par
$[x,y,z] := (x,y,z)-(-1)^{\bar{x}\bar{y}}(y,x,z)$,\\
one checks that $[x,y,z]$ expresses as (5.1) and  $(A, [ \cdot, \cdot, \cdot ], {\alpha}^2)$ turns out to be a Hom-Lie supertriple 
system. \hfill $\square$ \\ 
\par
We can now prove the main result of this section. \\
\par
{\bf Theorem 5.1} {\it The supercommutator Hom-algebra of any multiplicative right Hom-alternative superalgebra is a Hom-Bol superalgebra}.\\
\par
{\it Proof} Let ${\mathcal A} := (A,*, \alpha)$ be a multiplicative right Hom-alternative superalgebra. Then, by Lemma 5.1, ${\mathcal A}^{+}$ is a 
Hom-Jordan superalgebra and Lemma 5.2 says that $(A, [ \cdot, \cdot, \cdot ], {\alpha}^2)$ is a Hom-Lie supertriple system, where $[ \cdot, \cdot, \cdot ]$
is defined by (5.1). Now define on $A$ the ternary product
\par
$\{x,y,z\} := (-1)^{\bar{x}(\bar{y} + \bar{z})} as^{J}_{\alpha}(y,z,x)$. \\
Then, using properties of right Hom-alternative superalgebras, we get that \\ $(A,[\cdot,\cdot],\{ \cdot,\cdot,\cdot \}, {\alpha}^2)$ is a
Hom-Bol superalgebra. \hfill $\square$ \\ 

\vspace{0.5cm}
A. Nourou Issa \\ D\'epartement de Math\'ematiques, Universit\'e d'Abomey-Calavi, \\ 01 BP 4521 Cotonou 01, BENIN \\
email: woraniss@yahoo.fr
\end{document}